\documentclass[12pt,a4paper]{article}

\usepackage{latexsym}
\usepackage{amsfonts}
\usepackage{amssymb}
\usepackage{amsmath}
\usepackage{ifthen}
\usepackage{theorem}
\usepackage{enumerate}
\usepackage{graphics}

\newcommand{\normalmargins}{
  \setlength{\oddsidemargin}{-0.4mm}
  \setlength{\evensidemargin}{-0.4mm}
  \setlength{\textwidth}{160mm}

  \setlength{\topmargin}{-0.4mm}
  \setlength{\headheight}{4mm}
  \setlength{\headsep}{7mm}
  \setlength{\textheight}{225mm}
  \setlength{\footskip}{11mm}

  \addtolength{\topsep}{-\parskip}
  \addtolength{\partopsep}{-\parskip}

  \setlength{\headheight}{0mm}
  \setlength{\headsep}{0mm}
  \setlength{\textheight}{236mm}
}

\normalmargins

\theoremstyle{plain}
\newtheorem{blank}{}
\newtheorem{theorem}{Theorem}[section]
\newtheorem{conjecture}[theorem]{Conjecture}
\newtheorem{corollary}[theorem]{Corollary}
\newtheorem{lemma}[theorem]{Lemma}

{\theorembodyfont{\rmfamily} 
  \newtheorem{definition}[theorem]{Definition}
  \newtheorem{definitions}[theorem]{Definitions}
  
  \newtheorem{remark}[blank]{Remark}
    \newtheorem{example}[theorem]{Example}
    \newtheorem{example*}[blank]{Example}
}

\newcommand{\ifnull}[3]{\def\nullstring{}\def\teststring{#1}\ifx\nullstring\teststring#2\else#3\fi}
\newenvironment{proof}[1][]{\begin{trivlist}\item\ifnull{#1}{\noindent{\sc Proof:\ }}{\noindent{\sc Proof #1:\ }}}{\hspace*{\fill}$\Box$\end{trivlist}}

\newcommand{\figscale}{1}
\newcommand{\figref}[1]{Figure~\ref{fig:#1}}

\newcommand{\fig}[2]{
  \begin{figure}[htb]\centering
    \ifnull{#1}{}{\scalebox{\figscale}{\includegraphics{#1.eps}}}
    \ifnull{#2}{}{\caption{#2}}
    \label{fig:#1}
  \end{figure}
}

\renewcommand{\leq}{\leqslant}
\renewcommand{\geq}{\geqslant}

\renewcommand{\succeq}{\succcurlyeq}

\newcommand{\p}{^{\prime}}
\newcommand{\pp}{^{\prime\prime}}

\newcommand{\inr}[1]{{#1}^{\circ}}

\newcommand{\fto}{\rightarrow}
\newcommand{\embed}{\hookrightarrow}
\newcommand{\htop}{h_{\mathrm{top}}}

\newcommand{\olz}{{\overline{0}}}
\newcommand{\sfrac}[2]{{#1/#2}}

\newcommand{\zo}{\mbox{}^0_1}
\newcommand{\N}{\mathbb{N}}
\newcommand{\Z}{\mathbb{Z}}
\newcommand{\Q}{\mathbb{Q}}
\newcommand{\R}{\mathbb{R}}
\newcommand{\restrict}[1]{|_{#1}}
\newcommand{\ind}{\mathrm{Ind}}

\begin{document}

\title{Forcing relations for homoclinic and periodic orbits of the Smale horseshoe map
  \footnote{Mathematics subject classification: 
     Primary:   37E30, % Homeomorphisms and diffeomorphisms of planes and surfaces (Low-dimensional dynamical systems)
     Secondary: % 37B10, % Symbolic dynamics (Topological dynamics)
                37C27, % Homoclinic and heteroclinic orbits (Smooth dynamical systems)
                37E25, % Maps of trees and graphs (Low-dimensional dynamical systems)
    }
}
\author{Pieter Collins
  \thanks{This work was partially funded by Leverhulme Special Research Fellowship SRF/4/9900172. 
          The author would like to thank Toby Hall for numerous comments and suggestions which were invaluable in the writing of this paper.}
%\and Toby Hall
%\\Department of Mathematical Sciences\\University of Liverpool\\Liverpool L69 7ZL, U.K.
}
\date{Draft \today}

\maketitle

\begin{abstract}
An important problem in the dynamics of surface homeomorphisms is determining the forcing relation between orbits.  
The forcing relation between periodic orbits can be computed using standard algorithms, though this does not give much information on the structure of the forcing relation. 
Here we consider forcing relations between homoclinic orbits, and their relationships with periodic orbits.  
We outline a general procedure for computing the forcing relation, and apply this to compute the equivalence and forcing relations for homoclinic orbits of the Smale horseshoe map.
We also outline a method for determining forcing relations from graph maps, which allows us to compute the forcing relation between the so-called star homoclinic orbits.
\end{abstract}

%*****************************************************************************************************************************************************************************

\section{Introduction}
\label{sec:introduction}

We consider the problem of computing the braid equivalance and forcing relations for homoclinic and heteroclinic orbits of surface diffeomorphisms, with emphasis on the orbits of the Smale horseshoe map.
We outline a general method for computing the equivalence and forcing relations for homoclinic and heteroclinic orbits, and show how to obtain the periodic orbits forced by a given homoclinic orbit.
The approach is based on the trellis theory developed in \cite{Collins2002INTJBC,CollinsPPDYNSYS}, which is in turn based on the approach to the computation of the forcing relation for periodic orbits given in \cite{BestvinaHandel1995TOPOL,FranksMisiurewicz1993}.
We apply these methods to compute properties of the forcing relation for orbits of the Smale horseshoe map.

The theory requires fluency with a wide variety of concepts.
Orbits of the Smale horseshoe map are described combinatorially in terms of their natural symbolic coding.
Homoclinic orbits have a naturally associated trellis, which can be considered to be a subset of a homoclinic tangle, and provides a powerful framework for studying problems of forcing and braid equivalence.
Each trellis has a natural one-dimensional representative map, closely related to the train tracks for periodic orbits \cite{CassonBleiler88} , and these give another convenient combinatorial representation and computational tool.
Finally, by thickening graphs to obtain thick graphs, we obtain a link back to homoclinic and periodic orbits as described by braid types.

Much of this work was motivated by the \emph{decoration conjecture} of de Carvalho and Hall \cite{deCarvalhoHall2002EXPMAT} concerning the forcing relation between periodic orbits of the Smale horseshoe map.
The conjecture postulates that periodic orbits are partitioned into ordered families which are linearly ordered by the forcing relation, and the forcing relation can be determined from a knowledge of the forcing relation between families.
At the top of each family is a limiting homoclinic orbit, and the ordering between the families is determined the forcing relation between the homoclinic orbits.
If the conjecture is true, being able to compute the forcing relation between families is therefore an important problem in determining the forcing relation.

The paper is organised as follows. 
In Section~\ref{sec:orbits}, we describe the combinatorics of periodic and homoclinic orbits of the Smale horseshoe map, and state the decoration conjecture.
In Section~\ref{sec:trellis} we introduce the results of trellis theory which are necessary for the calculations.
In Section~\ref{sec:equivalenceforcing} we describe the algorithm used to compute the braid equivalence relation and the forcing relations and give a complete description of the forcing relation for homoclinic orbits of the Smale horseshoe map with short cores.
Finally, in Section~\ref{sec:stargraph}, we show how to compute the forcing relation for homoclinic orbits given in terms of their graph representatives, from which we deduce the forcing relation on the \emph{star homoclinic orbits} using results of \cite{deCarvalhoHallPPTOPOL}.

\enlargethispage{\baselineskip}

%*****************************************************************************************************************************************************************************

\section{Horseshoe Orbits}
\label{sec:orbits}

\newcommand{\rot}{\circlearrowleft}

%-----------------------------------------------------------------------------------------------------------------------------------------------------------------------------

\fig{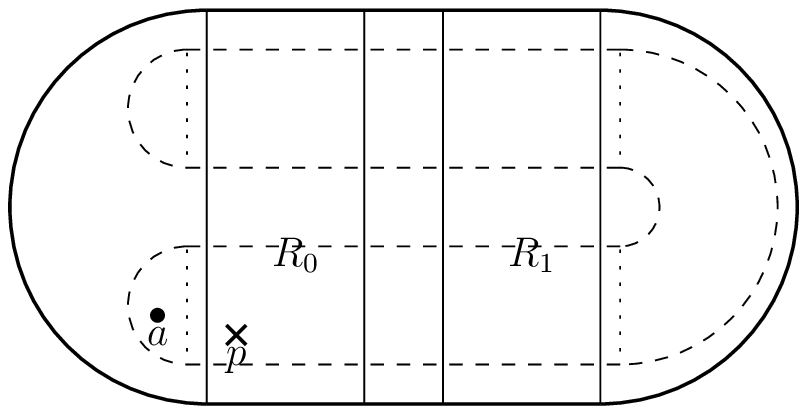}{The Smale horseshoe map}
In this paper, the model of the Smale horseshoe map $F\colon S^2\to
S^2$ depicted in~\figref{smalehs} will be used. The stadium-shaped
domain shown, consisting of two hemispheres and a rectangle $R$, 
is mapped into itself as in an orientation-preserving way as 
indicated by the dotted lines, with a stable fixed point $a$ in the 
left hemisphere. The map is then extended to a 
homeomorphism of $S^2$ with a repelling fixed point at $\infty$ whose
basin includes the complement of the stadium domain. 
The saddle fixed point of negative index (i.e. positive eigenvalues) is denoted~$p$.

The non-wandering set $\Omega(F)$ consists of the fixed points $a$ and
$\infty$, together with a Cantor set
\[\Lambda=\{x\in S^2\,:\,F^n(x)\in R \text{ for all }n\in\Z\}.\]

Since $\Lambda$ is contained in the union of the rectangles $R_0$ and
$R_1$, symbolic dynamics can be introduced in the usual way, providing
an \emph{itinerary} homeomorphism
\[k\colon \Lambda\to\Sigma_2=\{0,1\}^\Z,\]
with the property that $\sigma(k(x))=k(F(x))$ for all $x\in\Lambda$
(where $\sigma\colon\Sigma_2\to\Sigma_2$ is the shift map).

The itinerary $k(x)$ of a point~$x$ is periodic of period~$n$ if and
only if~$x$ is a period~$n$ point of~$F$. The following definition
makes it possible to describe periodic \emph{orbits} uniquely:

\begin{definition}[Code]
The \emph{code} $c_P\in\{0,1\}^n$ of a period~$n$ orbit $P$ of $F$ is
given by the first $n$ symbols of the itinerary of the rightmost point
of~$P$.
\end{definition}

Since the unimodal order on $\Sigma_2^+=\{0,1\}^\N$ reflects the
normal horizontal ordering of points, the elements of $\{0,1\}^n$
which are codes of period~$n$ orbits are those which are \emph{maximal}
in the sense of the following definition:

\begin{definitions}[Unimodal order, maximal word]
The \emph{unimodal order} $\prec$ on $\Sigma_2^+$ is defined as
follows: if $s=s_0s_1\ldots,t=t_0t_1\ldots\in\Sigma_2^+$ have
$s_n\not=t_n$, but agree on all earlier symbols, then $s\prec t$ if
and only if $\sum_{i=0}^n s_i$ is even. A word $w\in\{0,1\}^n$ is {\em
maximal} if $\sigma^j(\overline{w})\prec \overline{w}$ for $1\le j<n$.
\end{definitions}

This paper is also concerned with orbits which are homoclinic to $p$,
and the term \emph{homoclinic orbit} will be used exclusively to mean
such orbits. The points of a homoclinic orbit therefore have
itineraries containing only finitely many $1$s, and can thus be
described as follows:
\begin{definition}[Core]
Let $H$ be a homoclinic orbit of the horseshoe. The \emph{core} of $H$
is the longest word in the itinerary of a point of $H$ which begins
and ends with $1$.
\end{definition}

The \emph{signature} of $H$ is equal to the length of the core minus one.
Thus, for example, the homoclinic orbit containing the point of itinerary $\overline{0}110\cdot0101\overline{0}$ has core $1100101$ and signature $6$.
The \emph{primary} homoclinic orbits are those with cores $1$ and $11$; these two orbits have the same homoclinic braid type,
 are forced by every other homoclinic orbit, but do not force any other periodic or homoclinic orbit. 
By contrast, the orbits with cores $111$ and $101$ will be shown to force all periodic and homoclinic orbits of the horseshoe (cf.~\cite{Handel1999TOPOL}).

%-----------------------------------------------------------------------------------------------------------------------------------------------------------------------------

The \emph{rational words} defined next will be of particular importance
in what follows.
\begin{definition}[Rational word]
Given a rational number $q=m/n\in(0,1/2]$ (with $(m,n)=1$), define the
\emph{rational word} $c_q\in\{0,1\}^{n+1}$ by

\[
  c_{q,i} = \left\{ \begin{array}{l} 1 \textrm{ if }
     \left(\frac{m(i-1)}{n},\frac{m(i+1)}{n}\right) \textrm{ contains
     an integer,} \\ 0 \textrm{ otherwise} \end{array} \right.
\]
for $0\le i\le n$.
\end{definition}

The word $c_{m/n}$ can be determined by drawing a line from
$(0,0)$ to $(n,m)$ in the plane.  There is a $1$ in the $i$th position
if the line crosses an integer value in $i-1<x<i+1$. Thus, for
example, $c_{3/10}=10011011001$ (\figref{code3_10}). Note in
particular that, with the exception of the initial and final symbols,
$1$s always occur in blocks of even length.

\fig{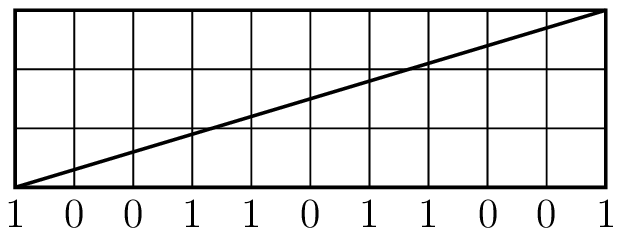}{The word $c_{3/10}$}

Using these rational words, it is possible to define the \emph{height}
of a horseshoe periodic orbit: this is a braid type invariant taking
rational values in $(0,1/2]$. The next lemma~\cite{Hall1994NONLIN} motivates the
definition: its proof is quite straightforward from the description
above, noting that lines of greater slope define words which are
smaller in the unimodal order.

\begin{lemma}
Let $q$ and $r$ be rationals in $(0,1/2)$ with $q<r$. Then $c_q0$,
$c_q1$, $c_r0$, and $c_r1$ are maximal, and
\[\overline{c_r0}\prec\overline{c_r1}\prec\overline{c_q0}\prec\overline{c_q1}.\]
\end{lemma}

\begin{definition}[Height]
Let $P$ be a horseshoe periodic orbit with code $c_P$. Then the {\em
height} $q(P)\in[0,1/2]$ of $P$ is given by
\[q(P)=\inf\{q\in\Q\cap(0,1/2]\,:\,q=1/2\text{ or }\overline{c_q0}\prec\overline{c_P}\}.\]
\end{definition}

That is, the itineraries $\overline{c_q0}$ are unimodally ordered
inversely to~$q$: the height of~$P$ describes the position in this
chain of the itinerary of the rightmost point of~$P$. Although it is
not obvious from the definition, $q(P)$ is always a strictly positive
rational, which can be computed algorithmically: see~\cite{Hall1994NONLIN} for
details.

There are several important classes of orbits which can be defined in
terms of the rational words:
\begin{itemize}
\item The \emph{rotation-compatible} periodic orbits of rotation number $m/n$
are period~$n$ orbits whose codes agree with $c_{m/n}$ up to the final
symbol. Thus, for example, the rotation-compatible orbits of rotation
number $3/10$ have codes $100110110\zo$. The two
rotation-compatible orbits of a given rotation number have the same
braid type (this is generally true for orbits whose codes differ only
in their final symbol), and force nothing but a fixed point.
\item The \emph{no bogus transition (NBT)} periodic orbits of height $m/n$ are the
period $n+2$ orbits whose codes start with $c_{m/n}$. Thus, for
example, the NBT orbits of height $3/10$ have codes
$10011011001\zo$. These orbits force all periodic orbits dictated by
the unimodal order, that is, all periodic orbits $P$ with
$\overline{c_P}\prec\overline{c_{m/n}\zo}$. In particular, by the
definition of height, they force all periodic orbits of height greater
than $m/n$.
\item The \emph{star} homoclinic orbits are those whose core is equal
to $c_{m/n}$ for some $m/n$. They are one of the main examples
considered in this paper, and will be studied in Section~\ref{sec:stargraph}.
\end{itemize}

A conjectural description~\cite{deCarvalhoHall2001JEURMS} of the structure of the forcing
relation on the set of homoclinic and periodic orbits of the horseshoe
can be given in terms of the \emph{decorations} of the orbits. The
definition depends on the following result~\cite{Hall1994NONLIN}:

\begin{lemma}
Let $P$ be a horseshoe periodic orbit of height $m/n$ which is not
rotation-compatible. Then $P$ has period at least $n+2$, and $c_P$
has $c_{m/n}$ as an initial word.
\end{lemma}

\begin{definition}[Decoration]
Let $P$ be a period $k$ horseshoe orbit of height $q=m/n$ which is not
rotation-compatible or NBT. Then the \emph{decoration} of $P$ is the
word $w\in\{0,1\}^{k-n-3}$ such that $c_P=c_q\zo w\zo$. (The empty
decoration is denoted $\cdot$). The decoration of a
rotation-compatible orbit is defined to be $\rot$, and that of an NBT
orbit to be $\ast$.

Let $H$ be a homoclinic horseshoe orbit whose core $c_H$ has length
$k\ge4$. Then the \emph{decoration} of $H$ is the word
$w\in\{0,1\}^{k-4}$ such that $c_H=1\zo w\zo 1$. The decoration of the
primary homoclinic orbits (with cores $1$ and $11$) is defined to be
$\rot$, and that of the homoclinic orbits with cores $111$ and $101$
is defined to be $\ast$.
\end{definition}

A periodic orbit of height $q$ and decoration $w$ is denoted $P_q^w$:
thus $c_{P_q^w}=c_q\zo w\zo$, provided that $w\not=\rot,\ast$. The
notation is justified by the result of~\cite{deCarvalhoHall2003NONLIN}, that all of the
(four or fewer) periodic orbits of height $q$ and decoration $w$ have
the same braid type. A homoclinic orbit of decoration $w$ is denoted
$P_0^w$: again, all the homoclinic orbits with the same decoration
have the same homoclinic braid type~\cite{deCarvalhoHall2002EXPMAT}.

Note that a decoration $w$ may not be compatible with all heights,
since $c_q\zo w\zo$ may not be a maximal word when~$q$ is large. 
For example, the word $c_q0 w0$ with $q=2/5$, $c_q=101101$ and $w=1001$
is not maximal, since the periodic sequence $\overline{101101010010}$ has
maximal word $100101011010$ with $q=1/3$, $c_q=1001$ and $w=101101$.

\begin{definition}[Scope]
The
\emph{scope} $q_w$ of a decoration $w$ is defined to be the supremum of
the values of $q$ for which a periodic orbit of height $q$ and
decoration $w$ exists. 
\end{definition}

The following result is from~\cite{deCarvalhoHall2002EXPMAT}:
\begin{lemma}
Let $w$ be a decoration. If $w=\rot$ or $w=\ast$, then
$q_w=1/2$. Otherwise, $q_w$ is the height of the periodic orbit
containing a point of itinerary $\overline{10w0}$.

If $0<q<q_w$, then each of the four words $c_q\zo w\zo$ (or each of the
two appropriate words whe $w=\rot$ or $w=\ast$) is the code of a
height $q$ periodic orbit, while if $1/2\ge q>q_w$, then none of the words is
the code of a height $q$ periodic orbit.
\end{lemma}

Thus, the set $\mathcal{D}_w$ of periodic and homoclinic
orbits of the horseshoe of decoration~$w$ is given by
\[\mathcal{D}_w=\{P_q^w\,:\,0\le q<=q_w\},\]
where the notation $<=$ indicates that $q=q_w$ is possible for some
decorations but not for others. Moreover, the union of the sets
$\mathcal{D}_w$ is the set of all periodic and homoclinic orbits of
the horseshoe.

The following states those parts of the decoration conjecture which
are relevant in this paper:

\begin{conjecture}[The Decoration Conjecture]
Given decorations $w$ and $w'$, write $w\succeq w'$ if the homoclinic
orbit $P_0^w$ forces the homoclinic orbit $P_0^{w'}$; and write $w\sim
w'$ if the two homoclinic orbits have the same homoclinic braid
type. Then
\begin{enumerate}[a)]
\item If $q<q'$ and $w\succeq w'$, then $P_q^w$ forces $P_{q'}^{w'}$.
\item $P_q^w$ and $P_{q'}^{w'}$ have the same braid type if and only
if $q=q'$ and $w\sim w'$.
\end{enumerate}
\end{conjecture}

In particular, a) implies that each family $\mathcal{D}_w$ of orbits
with a given decoration is linearly ordered by forcing, the order
being the reverse of the usual order on heights.

If this conjecture is true, then determining whether or not one
horseshoe periodic orbit forces another, and whether or not two
horseshoe periodic orbits have the same braid type, depends on being
able to carry out the corresponding computations for homoclinic
orbits. The main purpose of this paper is to describe a method using
which such computations on homoclinic orbits can be carried out. In
addition, the forcing relation on the set of star homoclinic orbits is
described completely.

%*****************************************************************************************************************************************************************************

\section{Horseshoe trellises}
\label{sec:trellis}

In this section those aspects of trellis theory which will be used
later are reviewed. Trellis theory is applicable in a much more
general setting (see~\cite{CollinsPPDYNSYS} for full details), but here the
key definitions and results are presented in a manner tailored for the
study of horseshoe trellises. All results stated in this section can
be found in~\cite{CollinsPPDYNSYS}.

The key ideas presented are as follows. A trellis is a finite portion
of the tangle of stable and unstable manifolds of a saddle fixed
point. Starting with the familiar tangle of the full horseshoe, the
full horseshoe trellis of signature~$n$ can be defined for each integer
$n\geq2$: it has longer and longer stable and unstable branches as $n$
increases. Given a horseshoe homoclinic orbit, the full horseshoe trellis
of appropriate signature can be \emph{pruned}, by removing as many
intersections as possible without disturbing the given homoclinic
orbit. This pruned trellis is a complete invariant of homoclinic braid
type, and so the technique can be used to determine whether or not two
given homoclinic orbits have the same braid type.

Given a trellis (and the action of a diffeomorphism on it), there is a
lower bound on the dynamics of any diffeomorphism which has such a
trellis. This minimal dynamics can be computed as the dynamics of a
tree map, using techniques similar to those of Bestvina and
Handel~\cite{BestvinaHandel1995TOPOL}. The dynamics forced by a given horseshoe homoclinic
orbit can thus be determined by finding the appropriate pruned
horseshoe trellis, and calculating the associated tree map.
By ``thickening'' the tree map, we obtain a canonical representative diffeomorphism which
 we show contains essentially all braid types forced by the homoclinic orbit.

%-----------------------------------------------------------------------------------------------------------------------------------------------------------------------------

\subsection{The full horseshoe trellis}

\begin{definition}[Trellis]
Let $f\colon S^2\to S^2$ be a diffeomorphism, and $p$ be a hyperbolic saddle fixed point of~$f$. 
Then a \emph{trellis for $f$} (\emph{at $p$}) is a pair $T=(T^U,T^S)$, where $T^U$ and $T^S$ are intervals in $W^U(f;p)$ and $W^S(f;p)$ respectively containing $p$.
(Here, $W^U(f;P)$ and $W^S(f;P)$ denote the unstable and stable manifolds, respectively, of $f$ at $p$.)
Given a trellis $T=(T^U,T^S)$, denote by $T^V$ the set of intersections of $T^U$ and $T^S$. 
The trellis is \emph{transverse} if all of its intersection points are transverse.
\end{definition}
Since all trellises considered in this paper will be transverse, and hence the word \emph{trellis} will be understood to mean \emph{transverse trellis}.

\begin{definition}[Segment]
Let $T=(T^U,T^S)$ be a trellis. A \emph{segment} of~$T$ is a closed
subinterval of either $T^U$ or $T^S$ with endpoints in $T^V$ but
interior disjoint from $T^V$. The segment is called \emph{unstable} or
\emph{stable} according as it is a subinterval of $T^U$ or of $T^S$.
\end{definition}

\begin{definition}[Region, bigon]
Let $T=(T^U,T^S)$ be a trellis. Then a \emph{region} of~$T$ is the closure 
of a component of $S^2\setminus(T^U\cup T^S)$.

A \emph{bigon} of~$T$ is a region bounded by two segments (one unstable
and one stable).
\end{definition}

Let $p$ be the fixed point of the horseshoe map~$F$ with code~$0$, and
let $W^U(F;p)$ and $W^S(F;p)$ be the unstable and stable manifolds
of~$p$. Let~$Q=(\ldots,q_{-2},q_{-1},q_0,q_1,q_2,\ldots)$ be the
homoclinic orbit of~$F$ with code~$\olz1\olz$, with the points labelled in
such a way that the itinerary $k(q_i)$ of $q_i$ is
$\sigma^i(\overline{0}\cdot 1\overline{0})$. Then $W^U(F;p)$ passes
successively through the points $\ldots,q_{-2},q_{-1},q_0,q_1,q_2,\ldots$, while $W^S(F;p)$
passes successively through the points $\ldots,q_{2},q_{1},q_0,q_{-1},q_{-2},\ldots$.

\begin{definition}[Full horseshoe trellis]
Given $i\ge0$ and $j\le 0$, denote by $\overline{T }^U_i$ an interval in
$W^U(F;p)$ with end intersections $p$ and $q_i$, and by $T^S_j$ the interval
in $W^S(F;p)$ with endpoints $p$ and $q_j$. 
For $n\ge0$, a \emph{full horseshoe trellis of signature~$n$} is a trellis $T=(T^U,T^S)$ for~$F$, 
 where $T^S=T^S_j$ for some $j$, and $T^U$ is a closed neighbourhood of $T^U_i$ in $W^U(F;p)$ for $i=j+n$ such that 
 all intersections of $T^U$ with $T^S$ lie $T^U_i$.
\end{definition}
It is clear that all full horseshoe trellises as defined above are differentiably conjugate;
 for definiteness, we will usually choose either $i=1$ or $i=\lfloor n/2\rfloor$.
We see that the endpoints of $T^S$ are intersection points, but the endpoints of $T^U$ are not.

\begin{example}
The full horseshoe trellis of signature~2 is depicted in ~\figref{smale}(a). 
The chaotics dynamics is supported in the regions labelled $R_0$ and $R_1$.
All points in $R_U$ are in the basin of the attracting fixed point $a$, and all points in the interior of $R_U$ are in the basin of the repelling point at infinity.
The point $r_0$ has itinerary $\overline{0}1\cdot01\overline{0}$, and the point $r_1$ has itinerary $\overline{0}1\cdot11\overline{0}$.
The full horseshoe trellis of signatures~3 is depicted in ~\figref{smale}(b). 
\fig{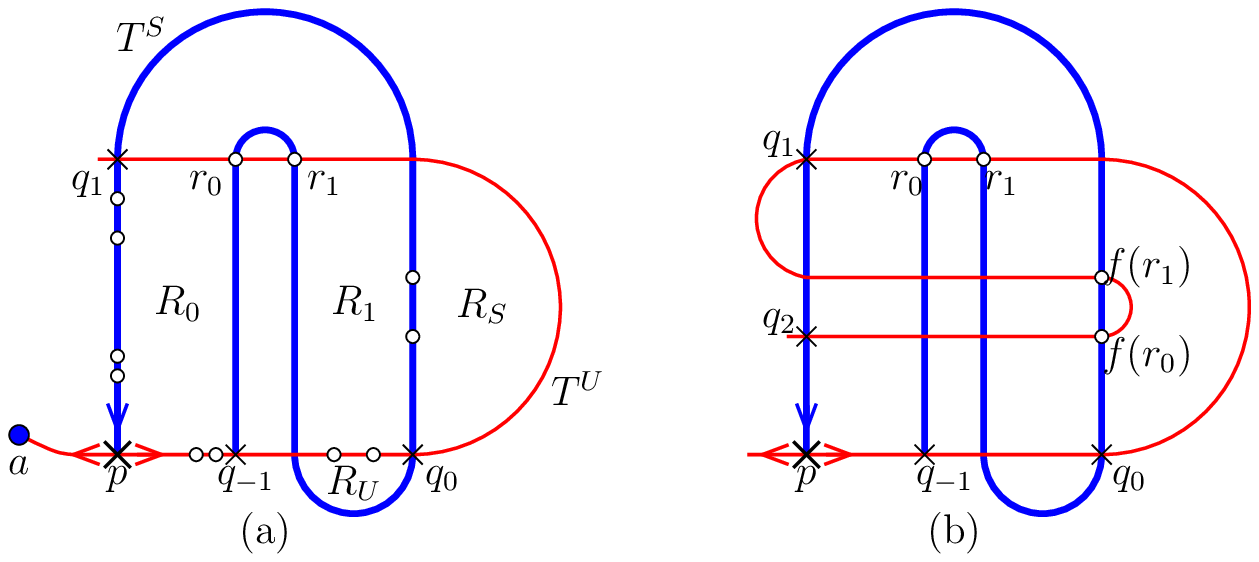}{(a) The full horseshoe trellis with signature $2$. (b) The full horseshoe trellis with signature $3$.}
\end{example}

The regions of a trellis can be used to introduce symbolic dynamics:
\begin{definition}[Itinerary]
Let $f$ be a diffeomorphism with trellis~$T$. Then a bi-infinite sequence $\ldots
R_{-2}R_{-1}R_0R_1R_2\ldots$ of regions of $T$ is an \emph{itinerary}
for an orbit $(x_i)$ of~$f$ if $x_i\in R_i$ for all~$i$.
\end{definition}

%-----------------------------------------------------------------------------------------------------------------------------------------------------------------------------

\subsection{Pruning isotopies and horseshoe trellises}
Given a trellis~$T$ for a diffeomorphism~$f$, a pruning isotopy is an
isotopy which removes the intersections on the boundary of one or more
bigons of $F$. To be more precise, it is an isotopy from $f$ to a
diffeomorphism~$f'$ which has a trellis $T'$ obtained from $T$ by
removing such intersections.
There are two possibilites; we can either remove both intersections of a single bigon, as depicted in~\figref{pruningisotopy}(a), or remove intersections from two neighbouring bigons, changing the orientation of the crossing at the remaining intersection, as depicted in ~\figref{pruningisotopy}(b).

\fig{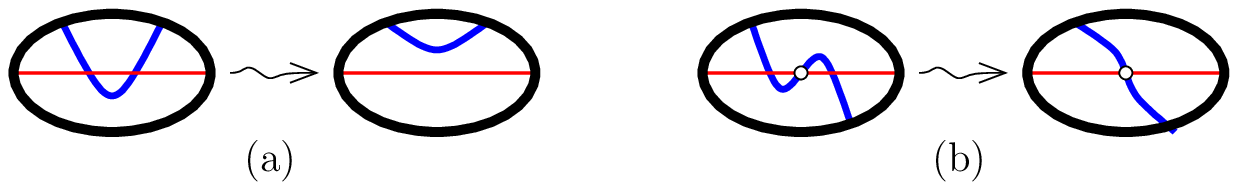}{The local effect of pruning isotopies on a trellis}
However, an isotopy of the diffeomorphism $f$ supported in some open set $U$ will also change the trellis outside of $U$.
If we are trying to reduce the number of intersections of $T$, we need to ensure that no other intersections are created when we remove intersections locally.
This gives rise to the notion of an \emph{inner bigon}
\begin{definition}[Inner bigon]
A bigon $B$ is \emph{inner} if $B\cap\bigcup_{n\in\Z}f^n(T^V)=B\cap T^V$.
i.e. A bigon is $B$ inner if the only intersections of $B$ with the orbits of the intersection points of the trellis are the vertices of $B$.
\end{definition}

The following result follows from the proof of Theorem~3.5 of \cite{CollinsPPDYNSYS}.
\begin{theorem}[Pruning away a bigon]
Let $T$ be a trellis of a diffeomorphism $f$.
Suppose either that $B$ is an inner bigon with vertices $v_0$ and $v_1$, or that $B_0$ and $B_1$ are inner bigons with a common vertex $v$ and other vertices $v_0$ and $v_1$ on different orbits.
Then there is a diffeomorphism $h:S^2\fto S^2$, which we can take to be supported on a neighbourhood $U$ of $B$ or $B_0\cup B_1$, such that
 $f\p=f\circ h$ has a trellis $T\p$ with the same intersections apart from those on the orbits of $v_0$ and $v_1$ under $f$.
\end{theorem}

Note that $f\p$ is isotopic to $f$, and has a trellis $T\p$ obtained by removing all the intersections of~$T$ contained in the orbit of $U$.

\begin{example}
\label{ex:prune}
\figref{smale3isotopy}~a) depicts the full horseshoe trellis of
signature~3, and a shaded neighbourhood~$U$ of a bigon~$B$, together
with its image. Pruning away the bigon~$B$ yields a diffeomorphism
$f'$ with trellis $T'$ as shown in~b). A further pruning isotopy
yields a diffeomorphism with the trellis depicted in~c).
\end{example}
\fig{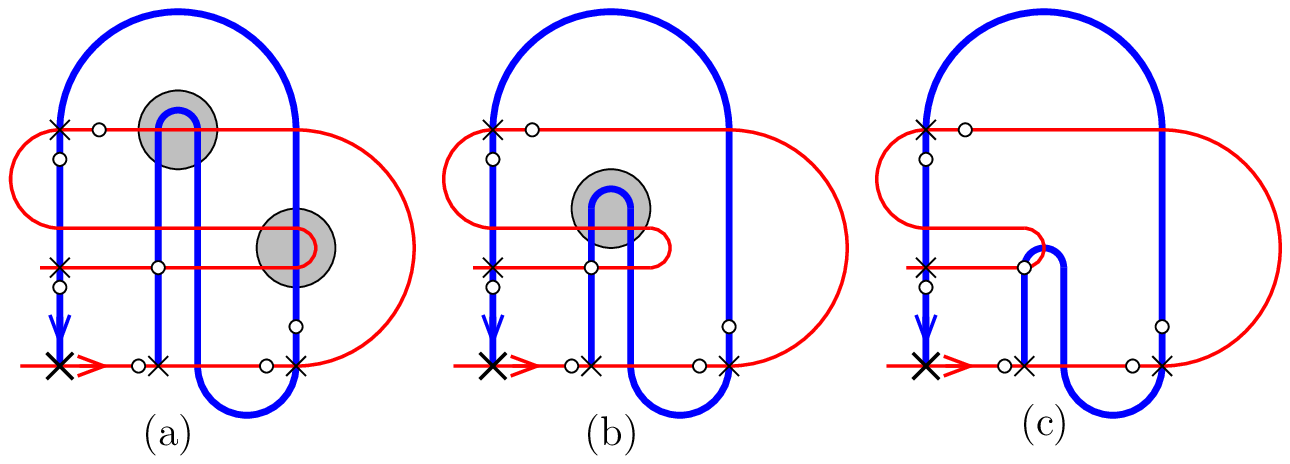}{Pruning away bigons in the horseshoe trellis}

A trellis obtained by pruning the full horseshoe trellis as in Example~\ref{ex:prune} is called a horseshoe trellis:
\begin{definition}[Horseshoe trellis]
A \emph{horseshoe trellis} is a trellis~$T$ obtained from the full horseshoe trellis by pruning away a sequence of bigons. 
\end{definition}
A horseshoe trellis can be associated to each homoclinic orbit of the horseshoe, by pruning away as many bigons as possible without touching the homoclinic orbit.
It is trivial that the signature of a horseshoe homoclinic orbit~$H$ is equal to the least integer~$n$
 such that $H$ is an intersection of the full horseshoe trellis of signature~$n$.

\begin{definition}[Trellis forced by a homoclinic orbit]
Let $H$ be a horseshoe homoclinic orbit of signature $n$.
The \emph{trellis forced by~$H$} is the trellis~$T$ obtained from the full horseshoe trellis of signature~$n$
 by pruning away as many bigons as possible which do not contain a point of~$H$.
\end{definition}

\begin{example}
The white circles in~\figref{smale3isotopy} represent points of the
homoclinic orbit~$H$ with code $\olz1001\olz$ (which thus has signature~$3$). The
trellis of \figref{smale3isotopy}c) is thus the trellis forced by this
homoclinic orbit. Note that every bigon has a point of~$H$ on its boundary.
\end{example}

This method makes it possible to determine whether or not two horseshoe homoclinic orbits have the same homoclinic braid type:
\begin{definition}[Trellis type]
Let $T$ and $T\p$ be horseshoe trellises for diffeomorphism $f$ and $f\p$ respectively.
We say that $(f;T)$ and $(f\p;T\p)$ have the same \emph{trellis type} if there is a diffeomorphism $g$ isotopic to $f$ relative to $T$, and a homeomorphism $h:S^2\fto S^2$ such that $h(T)=T\p$ and $h^{-1}\circ f\p\circ h=g$.
\end{definition}
We denote the trellis type containing $(f;T)$ by $[f;T]$.

For horseshoe trellises, the trellis type is determined by the geometry of the trellis:
\begin{theorem}
Let $T$ and $T\p$ be horseshoe trellises for diffeomorphism $f$ and $f\p$ respectively.
Then $(f;T)$ and $f\p;T\p)$ have the same trellis type if and only if $T$ and $T\p$ are diffeomorphic.
\end{theorem}
\begin{proof}
It suffices to consider the case $T=T\p$.
Since the points with itinerary $\overline{0}1\overline{0}$ all lie on a single homoclinic orbit,
 we can deduce the action of $f$ and $f\p$ on all the vertices of $T$ from the action of this orbit simply by counting vertices in each fundamental domain.
The result follows since all regions of $T$ are simply-connected, so the isotopy class is determined by the action on the segments.
\end{proof}
Since a horseshoe trellis type is fully determined by the geometry of the trellis, we define the type of a horseshoe trellis $T$ to be the type of $(f;T)$ for any diffeomorphism $f$ with trellis $T$ which can be obtained by pruning away bigons.

Horseshoe trellises $(f;T)$ and $(f\p;T\p)$ have the same type if and only if the trellises $T$ and $T\p$ are homeomorphic,
 and this occurs if and only if the orderings of the intersections on the stable and unstable manifolds are the same.
\begin{definition}
Let $(f;T)$ be a horseshoe trellis, with intersections $T^V=\{v_i:i=0\ldots n-1\}$
such that $v_i <_u v_{i+1}$ (i.e. $v_i$ is closer to $p$ along the unstable manifold).
Then the \emph{relative ordering} of the unstable and stable manifolds is the 
permutation $\pi_T$ such that $v_{\pi_T(i)} <_s v_{\pi_T(j)}$ if and only if $\pi_T(i)<\pi_T(j)$.
\end{definition}

The following result gives a computable criterion for the equivalence of horseshoe trellises.
\begin{theorem}
Let $T$ and $T\p$ be horseshoe trellises. Then $T$ and $T\p$ have the same trellis type if and only if $\pi_{T}=\pi_{T\p}$.
\end{theorem}

The following result is immediate from Theorem~3.5 of \cite{CollinsPPDYNSYS}, and shows that a homoclinic braid type is determined by the geometry of the trellis obtained by pruning up to the given orbit.
In particular, homoclinic orbits can only have the same braid type if they have the same signature (i.e. their cores have the same length).
\begin{theorem}
\label{thm:equiv}
Let $H$ and $H\p$ be horseshoe homoclinic orbits of signatures $n$ and $n\p$,  and let $T$ and $T\p$ be the trellises of signature $m\geq\max\{n,n\p\}$ forced by them. 
Then $H$ and $H\p$ have the same homoclinic braid type if and only if $T$ and $T\p$ have the same trellis type.
\end{theorem}

%-----------------------------------------------------------------------------------------------------------------------------------------------------------------------------

\subsection{The dynamics forced by a trellis}
The reason for the terminology `trellis \emph{forced} by a homoclinic
orbit~$H\p$ is that any diffeomorphism having a homoclinic orbit of the
homoclinic braid type of~$H$ has a fixed point and an associated
trellis of the given trellis type. Moreover, it is straightforward to
compute the dynamics forced by a trellis map, using techniques similar
to those of Bestvina-Handel~\cite{BestvinaHandel1995TOPOL} to represent this forced
dynamics by a graph map. It follows that the dynamics of this graph
map is forced by the homoclinic orbit. These intuitive notions are
made precise in this section.

\begin{definition}[Compatible tree]
Let $T$ be a horseshoe trellis, and $G$ be a tree embedded in
$S^2$. Then $G$ is \emph{compatible} with $T$ if
\begin{enumerate}[a)]
\item $G$ is disjoint from $T^U$.
\item $G$ intersects each segment of $T^S$ exactly once. The vertices of
$G$ are disjoint from $T^S$, and each edge of~$G$ intersects~$T^S$ at
most once.
\end{enumerate}
\end{definition}
The edges containing points of~$W$ are called \emph{control edges}. 
\begin{definition}[Compatible tree map]
Let $(f;T)$ be a horseshoe trellis map, and let~$G$ be a tree
compatible with~$T$. Let $W=G\cap T^S$.
A tree map $g\colon(G,W)\to(G,W)$ is \emph{compatible with $(f;T)$} if
\begin{enumerate}[a)]
\item $g(w_i)=w_j$ whenever $w_i,w_j\in W$ are the intersections of
$G$ with stable segments $S_i,S_j$ satisfying $f(S_i)\subseteq S_j$.
\item $g$ maps each control edge to a control edge.
\end{enumerate}
\end{definition}
\begin{definition}[Tree representative]
Let $(f;T)$ be a horseshoe trellis map, and $g\colon(G,W)\to(G,W)$ be
a tree map compatible with $(f;T)$. Then~$g$ is the \emph{tree
representative} of $(f;T)$ if
\begin{enumerate}[a)]
\item Every valence~$1$ or~$2$ vertex of $G$ is the endpoint of a
control edge.
\item $g$ is locally injective away from control edges (i.e. every
$x\in G$ has a neighbourhood~$U$ in $G$ such that if $y_1,y_2$ are
distinct points of $U$ with $g(y_1)=g(y_2)$, then at least one of
$y_1$ and $y_2$ lies in a control edge).
\item $g$ is piecewise linear on each edge of $G$.
\end{enumerate}
\end{definition}

An algorithm for computing the tree representative of a horseshoe
trellis map (and, more generally, the tree representative of an
arbitrary trellis map) can be found in~\cite{Collins2002INTJBC}. The
reason for the use of control edges is technical, and primarily
concerned with the details of the algorithm. Since the set of control
edges is invariant under the tree representative, and since the
concern here is with non-wandering dynamics, the tree representative
can be simplified by collapsing control edges to points.

\begin{definition}[Topological tree representative]
Let $g\colon(G,W)\to(G,W)$ be the tree representative of a horseshoe
trellis map $(f;T)$. Then the \emph{restricted tree representative} of
$(f;T)$ is $g\colon(\tilde G,\tilde W)\to(\tilde G,\tilde W)$, where
$\tilde G=\bigcap_{n=0}^\infty g^n(G)$, and $\tilde W=\tilde G\cap
W$. The \emph{topological tree representative} is obtained by collapsing
all control edges to points.
\end{definition}

The restricted tree representative contains all of the non-wandering
points of the tree representative, and so carries all of its
topological entropy. The topological tree representative is essentially
unique:
\begin{theorem}
Let $H$ and $H\p$ be horseshoe homoclinic orbits, and let $g$ and $g\p$
be the associated topological tree representatives. Then $g$ and $g\p$ are
conjugate if and only if $H$ and $H\p$ have the same
homoclinic braid type.
\end{theorem}
The proof is a simple corollary of the fact that the graph representative of a trellis type is unique.

\begin{example}
\fig{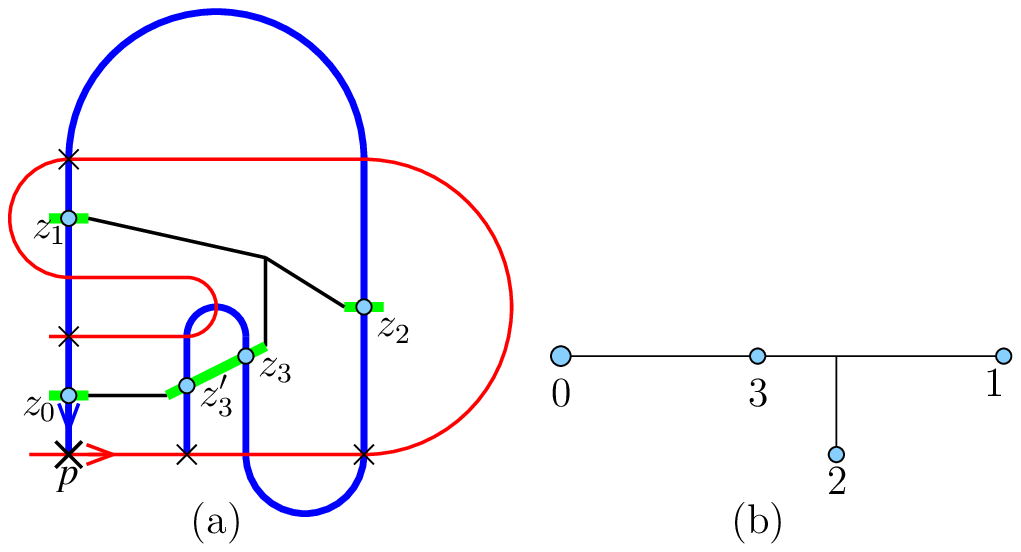}{(a) The restricted tree representative of the
trellis forced by $\overline{0}1111\overline{0}$. (b) The
topological tree representative.}  

Figure~\ref{fig:trellisgraph}(a) depicts the restricted tree representative of the trellis of signature $3$ forced by the homoclinic orbit $\overline{0}1111\overline{0}$.  
The control edges are labelled $z_0$, $z_1$, $z_2$, $z_3$ and $z_3\p$, and map under $g$ as
  \[ z_3, z_3\p \mapsto z_2 \mapsto z_1 \mapsto z_0 \mapsto z_0.  \] 
The topological tree representative of the orbit is shown in Figure~\ref{fig:trellisgraph}(b).  
\end{example}

The following conventions are used in labelling topological tree
representatives. The fixed point corresponding to the fixed point~$p$
of the horseshoe is labelled~$0$, and a preimage $x$ of~$p$ is
labelled with the least integer~$n$ satisfying $g^n(x)=p$. Thus, each
point labelled~$n$ is mapped to a point labelled~$n-1$. 
Where necessary, we use primes to distinguise $n$th preimages. 
We do not necessarily label all $n$th preimages of $v_0$, but will always
label points of the topological tree represesntative which are valence~$1$ vertices 
or at the \emph{fold vertices} at which the tree map is not locally injective.
In many cases, this labelling alone is enough information to determine the
tree map $g$.

\begin{example}
\fig{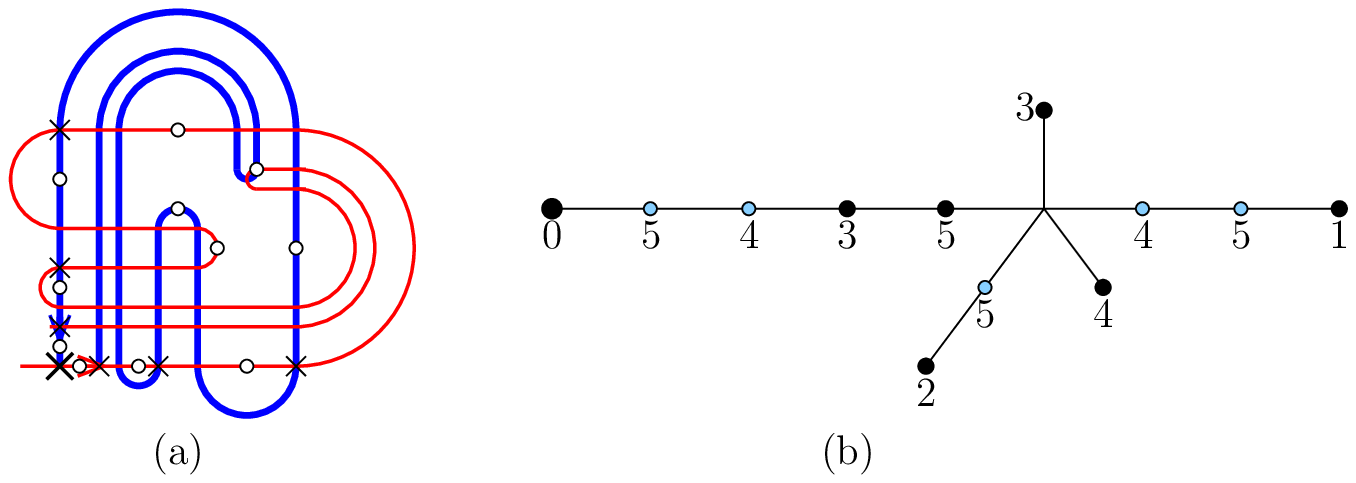}{Trellis forced by $\overline{0}111111\overline{0}$
and its topological tree representative.}  Figure~\ref{fig:trellis11}
shows the trellis forced by $\overline{0}1\zo11\zo1\overline{0}$ and
its topological tree representative.  The fold points and valence~1
vertices are marked with black dots.
\end{example}

The following theorem is Theorem~5.3 in \cite{CollinsPPDYNSYS}.
\begin{theorem}
Let $f$ be a diffeomorphism with a horseshoe trellis $T$ which has restricted tree representative~$g$. 
Then
\begin{enumerate}[a)]
\item For every orbit of~$g$ there is an orbit of~$f$ with the same itinerary.
\item For each periodic orbit of~$g$, there exists a periodic orbit of~$f$ with the same period and itinerary.
\item If $Y$ is an orbit of~$g$ which is homoclinic to the fixed control edge $z_0$ then there is a homoclinic orbit~$X$ of $f$ with the same itinerary as $Y$.
\item $h_{\text{top}}(f)\ge h_{\text{top}}(g)$.
\end{enumerate}
\end{theorem}
These results also hold for the topological tree representative, since the orbits of the restricted tree representative project to the topological tree representative.

The previous theorem shows that the topological tree representative of a trellis type give a good description of the orbits up to itinerary,
 but we often also want information about the braid type of the orbits, which is not directly given by the tree representative.
Following Franks and Misiurewicz \cite{FranksMisiurewicz1993} we ``thicken'' the topological tree representative to obtain the \emph{thick tree representative}.
The \emph{thick graph} $\widehat{G}\subset \R^2$ is a set consisting of \emph{thick vertices} $\widehat{V}$ and \emph{thick edges} $\widehat{E}$.
There is a bijection $\widehat{\cdot}$ between vertices and edges of $G$ and thick vertices and thick edges of $\widehat{G}$, a projection $\pi:\widehat{G}\fto G$ taking $\widehat{v}$ to $v$ and $\widehat{e}$ to $e$ which preserves the stable leaves.
There is also an embedding $i:G\embed\widehat{G}$ such that $i(v)$ lies in the interior of $\widehat{v}$, and $i(e)$ intersects $\widehat{e}$ in an unstable leaf, and intersects no other thick edges of $\widehat{G}$.

A \emph{thick tree representative} $\widehat{g}$ of $g$ is an embedding $g:\R^2\fto \R^2$ such that $\widehat{G}$ is a global attractor of $\widehat{g}$.
$\widehat{g}$ is a contraction on the set of thick vertices, uniformly contracts leaves of the stable foliation of the set of thick edges $\widehat{E}$ and uniformly expands the unstable foliation of $\widehat{E}$.
Further, $\widehat{g}(\widehat{v})\subset \widehat{g(v)}$ for each vertex $v$ of $G$, and $\widehat{g}(\widehat{e})$ crosses the same thick edges of $\widehat{G}$ in the same order as $g(e)$ crosses edges of $g$.
The thick tree representative is unique up to topological conjugacy.

Not all maps defined on trees embedded in $S^2$ give rise to thick tree maps, as edges may cross each other, resulting in a lack of embedding.
However, it is always possible to thicken the topological tree representative of a trellis type, and hence obtain the \emph{thick tree representative}.
It can be shown, though we do not do so here, that a horseshoe trellis type $T$ with a single homoclinic forcing orbit has is a trellis for its thick tree representative $\hat{g}$.
In general the thick tree representative may not have a trellis of the same type as the original diffeomorphism. 

We now show that almost every periodic orbit of $\widehat{g}$ may be continued to a periodic orbit of $f$.
This allows us to go full circle from homoclinic orbits through trellises, tree representaves, topological tree representatives and thick tree representatives and finally back to the original homoclinic orbit.
We say a periodic orbit of $\widehat{g}$ with least period $n$ is \emph{essential} if it is the only orbit in its Nielsen class.
It is easy to see that all periodic orbits of $\widehat{g}$ are essential except for those which are in the same Nielsen class as an attracting periodic orbit of $\widehat{g}$.
The only attracting periodic orbits of $\widehat{g}$ are contained in the vertices of $\widehat{G}$, and each attracting periodic point has one periodic point in the same Nielsen class in each adjoining thick edge.

\begin{theorem}
For each collection of essential periodic orbits of $\widehat{g}$, there is a collection of periodic orbits of $f$ with the same braid type.
\end{theorem}
\begin{proof}
\fig{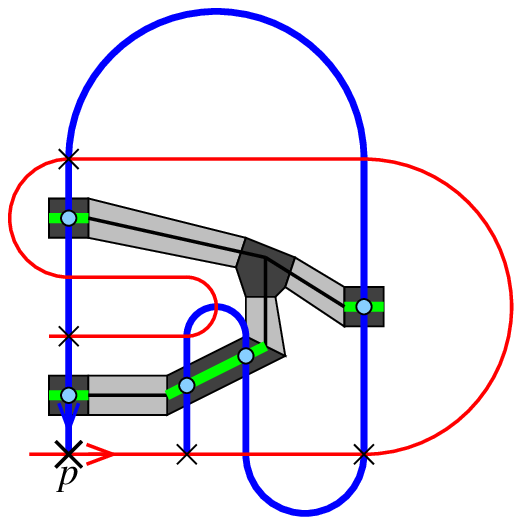}{A trellis $T$ with its thick graph representative $\widehat{G}$.}
Let $(g;G,W)$ be the tree representative of $[f;T]$.
Let $\widehat{G}$ be a compact neighbourhood of $G$ which deformation-retracts onto $G$, and $U$ a neighbourhood of $T^U$ which deformation-retracts to $T^U$ relative to $T^S$.
Then $K=R^2\setminus U$ deformation-retracts onto $G$, and further, there is an isotopy $h_t$ such that $h_t(T^U)=T^U$, $h_t(T^S)=T^S$ and $h_t(K)\subset \inr{\widehat{G}}$.
We can take $\widehat{G}$ to be a neighbourhood of $G$ with a thick-tree structure and $K$ to contain $f(\widehat{G})$.
Then $h_1\circ f$ maps $\widehat{G}$ into its interior and preserves the stable leaves.
We can then isotope in $\widehat{G}$ to obtain a diffeomorphism $\widehat{f}$ preserving the thick-tree structure, so that $\widehat{f}\restrict{\widehat{G}}=\widehat{g}$.

Since $h_t$ only preserves $T^S$ and $T^U$ set-wise, the sets $T^U$ and $T^S$ are still invariant for $h_1\circ f$, but will not be the stable and unstable manifolds of $T^P$.
However, $T^S$ is still invariant for $h_t\circ f$, so $\widehat{f}$ is isotopic to $f$ through diffeomorphisms $f_t$ for which $f_t(T^U)\subset T^U$ and $f_t(T^S)\subset T^S$.
This is sufficient to ensure that the essential Nielsen classes of $f$ and $\widehat{f}$ can continued through the isotopy.
Every periodic point of the topological tree representative $g$ of $[f;T]$ is the only periodic point in an essential Nielsen class.
Since $g$ is exact homotopy equivalent to $\widehat{f}$, these points lift to essential Nielsen classes of $\widehat{f}$ with the same Nielsen number.
For each periodic point $y$ of $g$, there is an essential Nielsen class $N$ related to $y$ by the projection onto $G$.
Further, $N$ contains a periodic point $\widehat{y}$ in $\widehat{G}$.
Since $\ind(N;\widehat{f})=\ind(y;g)=\ind(\widehat{y};\widehat{G})$, the point $\widehat{y}$ is in an essential Nielsen class of $\widehat{f}$ as well as of $\widehat{g}$.
Hence the periodic orbit $\widehat{y}_i$ of $\widehat{y}$ under $\widehat{f}$ can be continued by isotopy to a periodic orbit $(x_i)$ of $f$ of the same period and braid type.
\end{proof}

\begin{remark}
We note that if $[f;T]$ the trellis type forced by a single horseshoe homoclinic orbit, then the thick tree representative $\widehat{g}$ has a trellis $\widehat{T}$ such that $[\widehat{g};\widehat{T}]=[f;T]$.
This follows since by the main theorem of \cite{CollinsPPMODEL}, we can find a diffeomorphism $\widehat{f}\in [f;T]$ such that $\widehat{f}$ has the same entropy as $f$.
The construction of this diffeomorphism shows that $\widehat{f}$ is conjugate to the thick tree representative.
The result need not be true for a trellis forced by a collection of homoclinic orbits.
\end{remark}

%*****************************************************************************************************************************************************************************

\section{Orbit equivalence and forcing}
\label{sec:equivalenceforcing}

In this section the techniques described in Section~\ref{sec:trellis}
are used to compute, for horseshoe homoclinic orbits with short cores,
the equivalence classes under the relation of having the same homoclinic braid
type, and the forcing relation.

%-----------------------------------------------------------------------------------------------------------------------------------------------------------------------------

\subsection{Horseshoe homoclinic orbits with the same homoclinic braid
type}
\label{sec:samehbt}

We have computed the trellis types forced by all horseshoe homoclinic orbits of
signature~12 or less, and applied Theorem~\ref{thm:equiv}
to determine which pairs have the same homoclinic braid
type. Table~\ref{tab:equivalences} presents the results for signatures
up to~9.

For orbits of signature up to $4$, the only equivalences are trivial;
two homoclinic orbits have the same homoclinic braid type if and only
if they have the same decoration.  However the orbits of signature $5$
with decorations $01$ and $10$, and codes $\olz1\zo01\zo1\olz$ and
$\olz1\zo10\zo1\olz$ respectively, are equivalent (see
\figref{trellis01}).
\fig{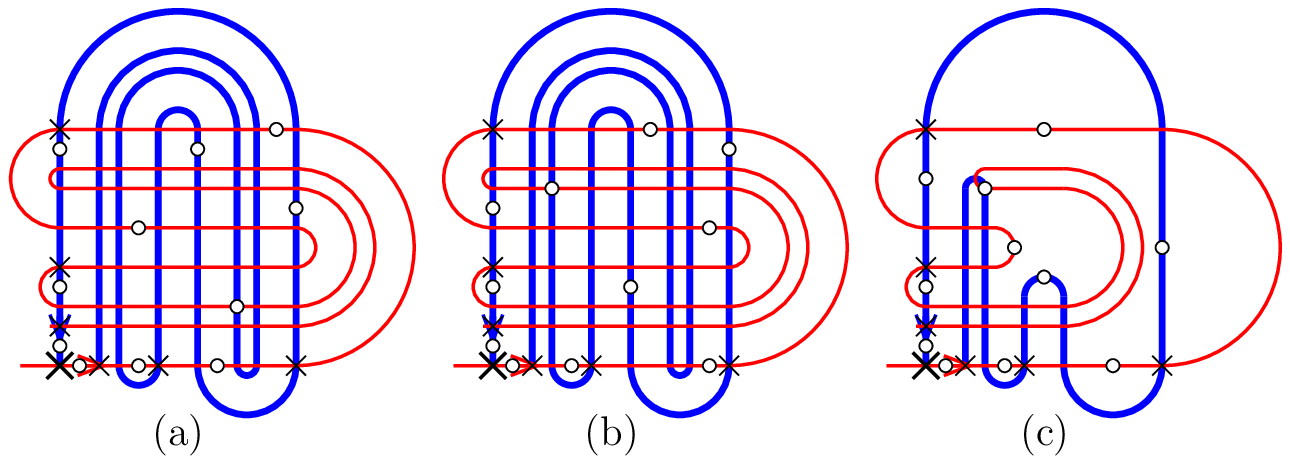}{The horseshoe trellis and homoclinic orbits with codes (a) $\overline{0}110111\overline{0}$ and (b) $\overline{0}111011\overline{0}$, and (c) the trellis they force.}

For homoclinic orbits of signature at most $7$, all orbits have the
same homoclinic braid type as their time reversal.  However, for
orbits of signature $8$, there are two pairs of orbits whose
homoclinic braid type differs from that of their time reversal.  This
is a counterexample to the conjecture that horseshoe orbits which are
time-reversals have the same braid type.

\fig{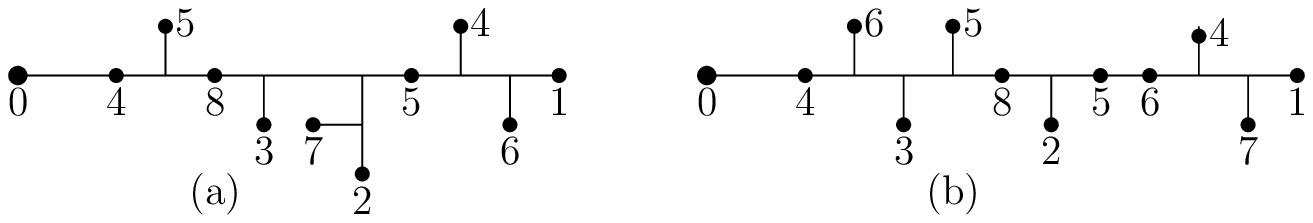}{Topological graph representatives for the homoclinic orbits with decoration (a) $01001$ and (b) $10010$.}

The trellises forced by the orbits with codes $\overline{0}1\zo01001\zo1\overline{0}$ and
$\overline{0}1\zo10010\zo1\overline{0}$ are not equivalent, even though these words are reverses
of each other.  The graph representatives are shown in
Figure~\ref{fig:graphreverse01001}.  The topological entropy of both
these orbits is $\log\lambda$, where $\lambda$ is the largest root of
the polynomial
\[
  \lambda^{13}-2\lambda^{12}+2\lambda^8+\lambda^7-4\lambda^5-2\lambda^4+2\lambda^2+2\lambda-2
  .
\]
Numerically, $\lambda_{\max}\approx1.845$, giving $\htop>0.612$.

\fig{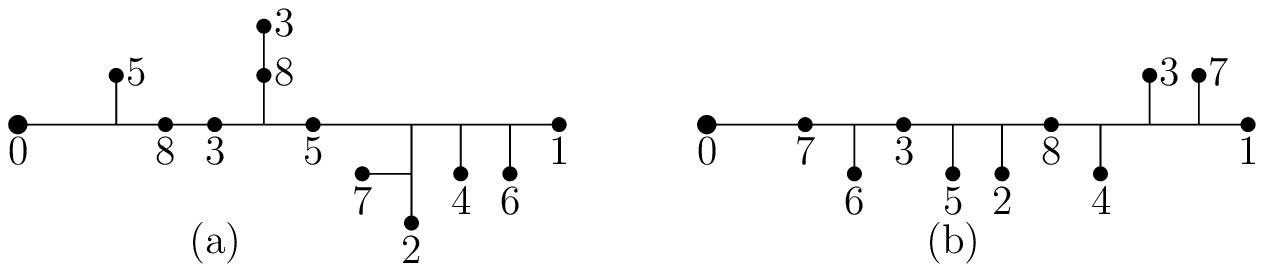}{Topological graph representatives for the homoclinic orbits with decoration (a) $11001$ and (b) $10011$.}

Another counterexample is given by the orbits with codes
$\overline{0}1\zo11001\zo1\overline{0}$ and $\overline{0}1\zo10011\zo1\overline{0}$.  The topological
entropy of both these orbits is $\log\lambda$, where $\lambda$ is the
largest root of the polynomial
\[ 
  \lambda^{13}-2\lambda^{12}+3\lambda^7-4\lambda^6+4\lambda^5+2\lambda^4+2\lambda-2
  .
\]
Numerically, $\lambda_{\max}\approx1.909$, giving $\htop>0.646$.

\begin{table}
\[
\begin{array}{|c|c|l|} \hline
  \textrm{Sig} & \textrm{Scope} & \textrm{Decorations} \\ \hline
   2 & \sfrac{1}{2} & * \\
  \hline
   3 & \sfrac{1}{3} & \cdot \\
  \hline
   4 & \sfrac{1}{4} & 0 \\
     & \sfrac{1}{2} & 1 \\
  \hline
   5 & \sfrac{1}{5} & 00 \\
     & \sfrac{1}{4} & 01,\ 10 \\
     & \sfrac{2}{5} & 11 \\
  \hline
   6 & \sfrac{1}{6} & 000 \\
     & \sfrac{1}{5} & 001,\ 100 \\
     & \sfrac{1}{3} & 011,\ 010,\ 110 \\
     & \sfrac{1}{2} & 111 \\
     & \sfrac{1}{2} & 101 \\
  \hline
   7 & \sfrac{1}{7} & 0000 \\
     & \sfrac{1}{5} & 0001,\ 1000 \\
     & \sfrac{1}{4} & 0011,\ 0010,\ 0100,\ 1100 \\
     & \sfrac{2}{7} & 0110 \\
     & \sfrac{1}{3} & 0111,\ 0101,\ 1110,\ 1010 \\
     & \sfrac{1}{3} & 1001 \\
     & \sfrac{2}{5} & 1101,\  1011 \\
     & \sfrac{3}{7} & 1111 \\
  \hline
   8 & \sfrac{1}{8} & 00000 \\
     & \sfrac{1}{6} & 00001,\ 10000 \\
     & \sfrac{1}{5} & 00011,\ 00010,\ 01000,\ 11000 \\
     & \sfrac{1}{4} & 00110,\  00100,\ 01100 \\
     & \sfrac{1}{4} & 00111,\ 00101,\ 11100,\ 10100 \\
     & \sfrac{1}{4} & 10001 \\
     & \sfrac{1}{3} & 01101,\ 01111,\ 01110,\ 01010,\ 01011,\ 11010,\ 11110,\ 10110 \\
     & \sfrac{1}{3} & 01001 \\
     & \sfrac{1}{3} & 11001 \\
     & \sfrac{1}{3} & 10010 \\
     & \sfrac{1}{3} & 10011 \\
     & \sfrac{3}{8} & 11011 \\
     & \sfrac{1}{2} & 11111 \\
     & \sfrac{1}{2} & 11101,\ 10111 \\
     & \sfrac{1}{2} & 10101 \\
   \hline
\end{array}
\]
\end{table}

\begin{table}
\[
\begin{array}{|c|c|l|} \hline
  \textrm{Sig} & \textrm{Scope} & \textrm{Decorations} \\ \hline
   9 & \sfrac{1}{9} & 000000 \\
     & \sfrac{1}{7} & 000001,\ 100000 \\
     & \sfrac{1}{6} & 000011,\ 000010,\ 010000,\ 110000 \\
     & \sfrac{1}{5} & 000110,\ 000100,\ 001000,\ 011000 \\
     & \sfrac{1}{5} & 000111,\ 000101,\ 111000,\ 101000 \\
     & \sfrac{1}{5} & 100001 \\
     & \sfrac{2}{9} & 001100 \\
     & \sfrac{1}{4} & 001101,\ 001110,\ 001010,\ 001001,\ 011100,\ 010100,\ 101100,\ 100100 \\
     & \sfrac{1}{4} & 001111,\ 001011,\ 110100,\ 111100 \\
     & \sfrac{1}{4} & 010001 \\
     & \sfrac{1}{4} & 110001 \\
     & \sfrac{2}{7} & 011001,\ 100110 \\
     & \sfrac{1}{3} & 011011,\ 011010,\ 011110,\ 010110,\ 110110 \\
     & \sfrac{1}{3} & 011111,\ 010111,\ 111110,\ 111010 \\
     & \sfrac{1}{3} & 011101,\ 010101,\ 101010,\ 101110 \\
     & \sfrac{1}{3} & 010010,\ 010011,\ 110010 \\
     & \sfrac{1}{3} & 110011 \\
     & \sfrac{1}{3} & 111001 \\
     & \sfrac{1}{3} & 101001 \\
     & \sfrac{1}{3} & 100101 \\
     & \sfrac{1}{3} & 100111 \\
     & \sfrac{2}{5} & 110111,\ 110101,\ 111011,\ 101011 \\
     & \sfrac{2}{5} & 101101 \\
     & \sfrac{3}{7} & 111101,\ 101111 \\
     & \sfrac{4}{9} & 111111 \\
  \hline
\end{array}
\]  
\caption{Equivalences for homoclinic orbits with signatures up to $9$.}
\label{tab:equivalences}
\end{table}

\subsection{The forcing relation}

The homoclinic orbits forced by a given homoclinic braid type can also
be computed using the methods of Section~\ref{sec:trellis}. Since the
trellis forced by a homoclinic orbit can be computed without
introducing new intersections, the (horseshoe) itineraries of
intersections which remain can be continued through the pruning. The
forced homoclinic braid types are precisely those which persist
through the pruning.

Note, though, that in order to show that the homoclinic braid type of
a homoclinic orbit~$H$ does {\em not} force that of a homoclinic
orbit~$H'$, it is necessary to show that {\em none} of the homoclinic
orbits of the same type as~$H'$ persist through the pruning. Thus it
is necessary to compute equivalences, as in Section~\ref{sec:samehbt},
in order to be able to compute the forcing relation.

\fig{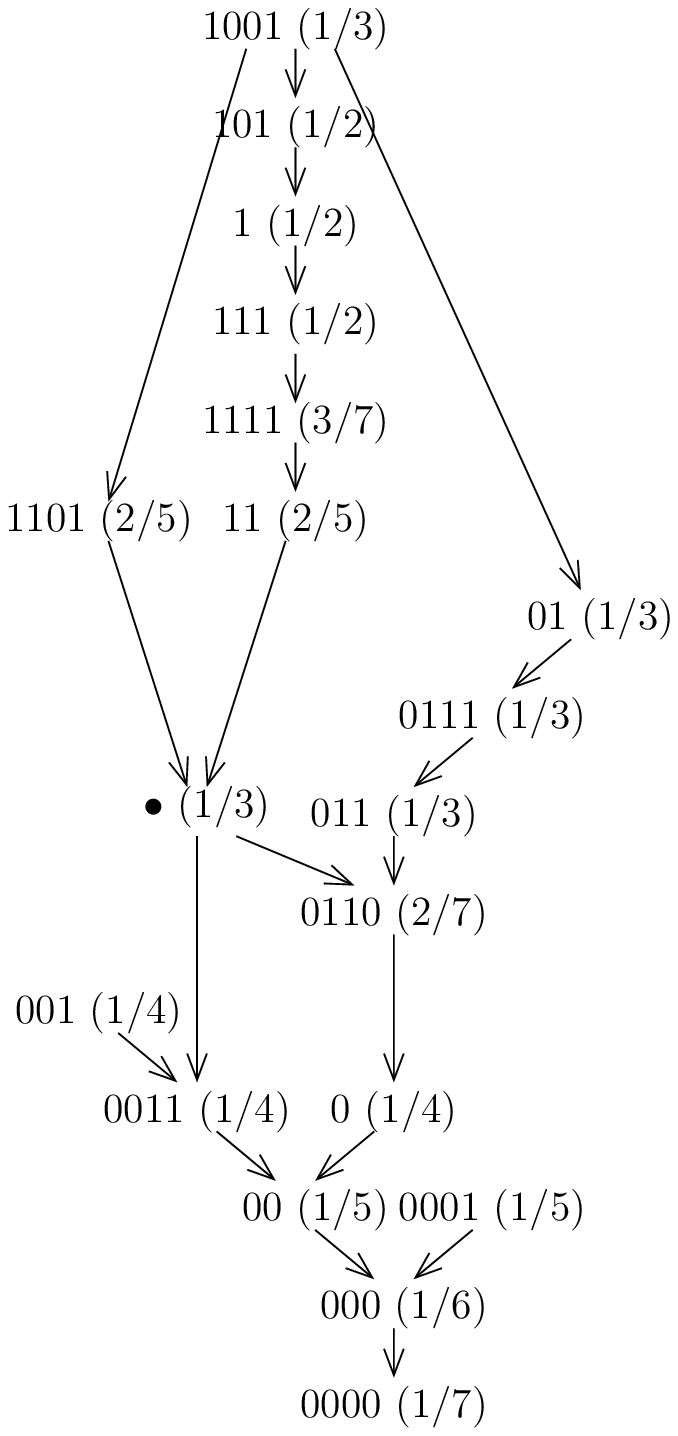}{Decorations and scopes of homoclinic braid types of
signature up to $7$ and the forcing relation between them.}

\figref{forcing} shows the forcing relation between horseshoe
homoclinic orbits of signature~7 or less. Each homoclinic orbit is
specified by its decoration, and the scope of the decoration is also
given. Only one decoration is given for each equivalence class of
homoclinic braid types (so, for example, the decoration $10$ is not
included, since it is equivalent to $01$).

In the remainder of the section, two examples are given to indicate
how \figref{forcing} has been computed.
\begin{example}\ \par
\fig{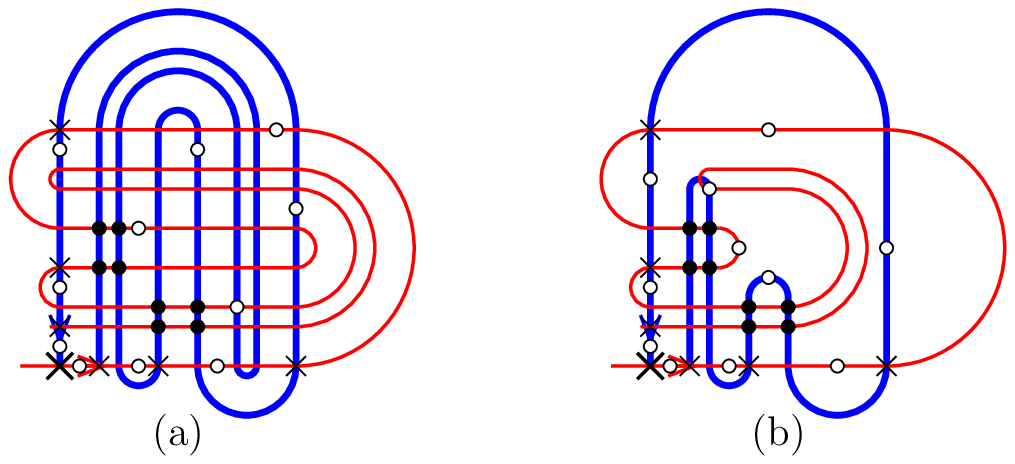}{The homoclinic orbit $\olz110111\olz$ forces the orbits $\olz1\zo0\zo1\olz$.}
Figure~\ref{fig:forcing01} shows the horseshoe trellis, on which the
orbit $P$ with code $\olz110111\olz$ and decoration $01$ has been marked with white dots,
and the intersections on the orbits with codes $\olz1\zo0\zo1\olz$
with black dots.  On performing a pruning isotopy to obtain the
trellis forced by the orbit $P$, it can be seen that the marked
intersections persist.  Therefore the homoclinic braid type with
decoration $01$ forces the homoclinic braid type with decoration $0$.
From Figure~\ref{fig:forcing01} it can also be seen that the only
other forced homoclinic braid type with signature less than $5$ is
that with decoration $00$.
\end{example}

\begin{example}\ \par
\fig{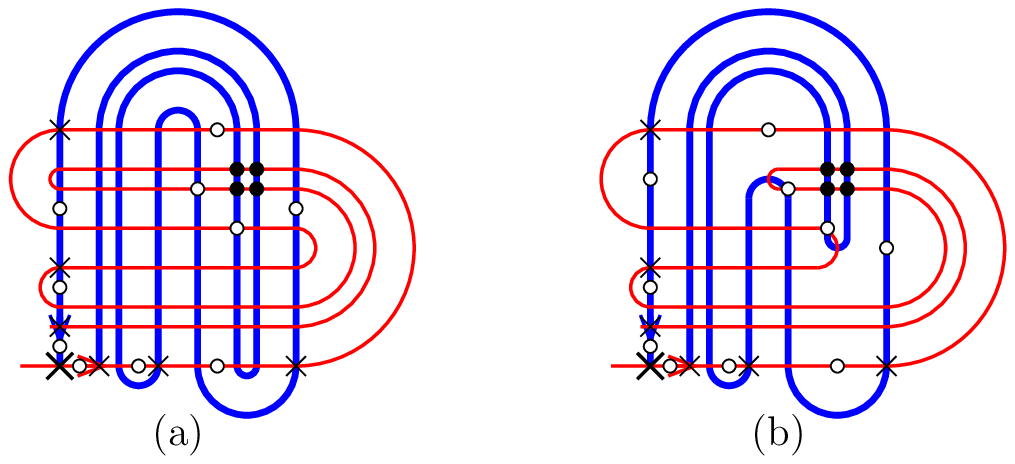}{The homoclinic orbit $\olz11111\olz$ forces the orbits
  $\olz1\zo11\zo1\olz$.}
\par
It is possible to compute braid types of arbitrary signature forced by any homoclinic orbit.
In Figure~\ref{fig:forcing1}, we show the trellis of signature $5$ forced by the homoclinic orbit $\olz11111\olz$ which has decoration $1$ and signature $4$, allowing us to compute the braid types of signature up to $5$ forced by the orbit.
In particular, we see that the orbit $\olz11111\olz$ forces the orbits $\olz1\zo11\zo1\olz$, so the homoclinic braid type with decoration $1$ forces the homoclinic braid type with decoration $11$.
The other homoclinic braid types with signature less than $5$ forced by the braid type with decoration $1$ have decorations $\bullet$, $0$ and $00$.
\end{example}

%*****************************************************************************************************************************************************************************

\section{Graph representatives and star orbits}
\label{sec:stargraph}

The methods described above can be used to compute the forcing relation between homoclinic braid types of arbitrarily high signature,
 but cannot be used to describe the forcing relation on infinite sets of homoclinic braid types.
Further, we would also like to be able to compute the forcing relation between periodic orbits.
In this section we use the topological tree representatives of trellis types to obtain more information on the forcing relations.
We need two operations, glueing and pulling tight, which can be defined both at the level of trees and thick trees.
The operations were used in \cite{deCarvalhoHall2001JEURMS,deCarvalhoHallPPTOPOL}, and are based on those of \cite{BestvinaHandel1995TOPOL,FranksMisiurewicz1993}.

%-----------------------------------------------------------------------------------------------------------------------------------------------------------------------------

\subsection{Forcing of graph representatives}
\label{sec:graphforcing}

We first show how to construct a forcing relation for topological tree representatives.
This allows us to perform calculations on the tree level, and only go to the level of homoclinic orbits and trellises to prove our forcing results.
We remind the reader that all tree maps considered in this section can be thickened to give thick tree maps.
\begin{definition}[Glueing and pulling tight]
A tree map $(g\p;G\p)$ is obtained from $(g;G)$ by \emph{glueing} if there is a surjective tree map $h$ such that $h\circ g=g\p\circ h$.
A tree map $(g\p;G)$ is obtained from $(g;G)$ by \emph{pulling tight} for every edge $e$ of $G$, $g\p(e)$ is a sub-path of $g(e)$.
\end{definition}
A tree map $(g\p;G\p)$ is \emph{forced} by a tree map $(g;G)$ if $(g\p;G\p)$ can be obtained from $(g;G)$ by glueing and pulling tight.
Note that if $(g\p;G\p)$ is forced by $(g;G)$ and $W\p$ is an invariant set of $g\p$, then $W=h^{-1}(W\p)$ is invariant under $g$.
The operations of glueing and and edge collapsing have a clear interpretation in terms of the thick tree representatives.

\begin{example}\ 
\fig{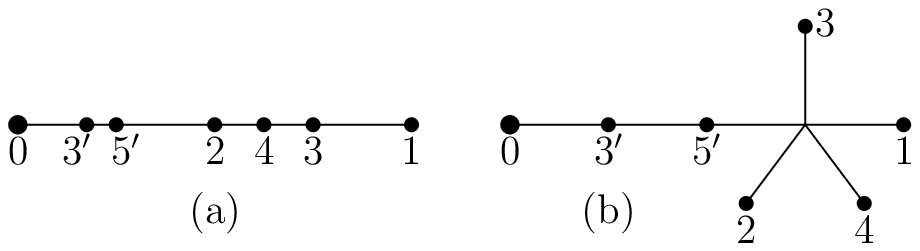}{(a) The tree map $g_{1/2}$ with an invariant set $W$. (b) The tree map $g_{2/5}$ with an invariant set $W_{2/5}$.}
The tree map $g_{1/2}$ depicted in \figref{graphforcing}(a) forces the tree map $g_{2/5}$ depicted in \figref{graphforcing}(b), since we can take $h$ to be a piecewise-linear map taking each point of the labelled invariant vertex set $W$ of $G_{1/2}$ to the point the invariant vertex set $W_{2/5}$ of $G_{2/5}$ with the same label.
\end{example}

\begin{theorem}
\label{thm:gluetighten}
Let $g$ be a tree map, and $g\p$ be a tree map which can be obtained from $g$ by the operations of glueing and pulling tight.
Let $f$ be the thick tree represenative of $g$, and $f\p$ the thick tree representative of $g\p$.
Then there is a homeomorphism $h$ such that $h\circ f$ is isotopic to $f\p\circ h$ relative to $h^{-1}(\Lambda\p)$, where $\Lambda\p$ is the nonwandering set of $g\p$.
\end{theorem}

\begin{proof}
Suppose $g\pp$ is obtained from $g$ by glueing, so there is a tree map $h:G\fto G\p$ such that $h\circ g=g\pp\circ h$.
Let $D$ and $D\p$ be the thick trees of $G$ and $G\p$ respectively, and let $\hat{h}:D\fto D\p$ be the map induced by $h$.
Let $f\pp$ be the thick tree representative of $g\pp$.
Then $f=\hat{h}^{-1}\circ f\pp \circ \hat{h}$ is a thick tree representative of $g$.

Now suppose $g\p$ is obtained from $g\pp$ by pulling tight.
Since we can realise pulling tight by a pruning isotopy, we can choose this to be constant on the nonwandering set $\Lambda\p$ of $f\p$.
Hence $\hat{h}\circ f\circ\hat{h}^{-1}$ is isotopic to $f\p$ relative to $\Lambda\p$, so $\hat{h}\circ f$ is isotopic to $f\p\circ\hat{h}$ relative to $\hat{h}^{-1}(\Lambda\p)$.
\end{proof}

An immediate corollary of this result shows how to obtain results on the trellis level.
\begin{corollary}
Suppose $(g;G,W)$ is the topological tree representative of a horseshoe trellis type $[f;T]$, where $f$ is conjugate to the thick tree representative of $g$.
Suppose further that $(g\p;G\p,W\p)$ is obtained from $(g;G,W)$ by glueing and pulling tight.
Then $(g\p;G\p,W\p)$ is the topological tree representative of a horseshoe trellis type $[f\p;T\p]$ which can be obtained from $[f;T]$ by pruning away bigons.
\end{corollary}

%-----------------------------------------------------------------------------------------------------------------------------------------------------------------------------

\subsection{The star homoclinic orbits}
\label{sec:star}

Of special interest are the star homoclinic orbits.
These are the horseshoe orbits for which the tree representative is the \emph{$n$-star}.
To describe the tree map $g_{m/n}$, we label the edges of $\Gamma_n$ by $e_i$, with edge $e_{i+m\bmod n}$ immediately following $e_i$ in the (anticlockwise) cyclic order at the centre vertex $v$.

\begin{remark}
The labelling described here is different from that of e.g. \cite{deCarvalhoHallPPTOPOL}, in which the edges are labelled $e_0,e_1,\ldots e_{n-1}$ cyclically.
We can think of this cyclic labelling as a \emph{geometric} labelling, as opposed to the \emph{dynamical} labelling used here.
\end{remark}

To define the map $g_{m/n}$ we let $r$ be (unique) integer with $1\leq r<n$ and $mr=-1\bmod n$.
(In other words, $r=-1/m\bmod n$. 
We define the map $g_{m/n}$ by
\[
  \begin{array}{rcl} 
    g_{m/n}(e_0) & = & e_0\bar{e}_{r}e_{r}\bar{e}_{2r}e_{2r}\cdots \bar{e}_2e_2 \\
    g_{m/n}(e_i) & = & e_{i-1} \textrm{ for } 1\leq i\leq n-1
  \end{array}
\]
The end of the edge $e_i$ is a valence-$1$ vertex $i$.
There are $m$ fold points, all of which are on $e_0$, and which have labels $(ir+1\mod n)\p$ for $1\leq i<n$.
The $m$th fold point is an inner fold with label $m\p$.
By convention we shall draw star trees so that the vertices $0$ and $1$ lie on a horizontal line.

\begin{example}
\fig{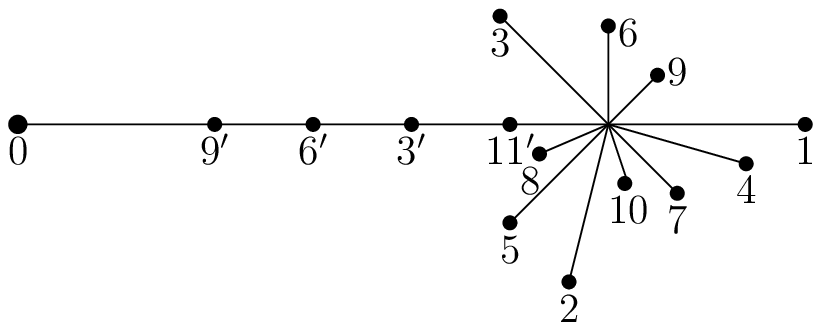}{The star map $g_{4/11}$.}
Consider the tree map $g_{4/11}$.
Then $r=-1/4\bmod 11=8$.
The edges have labels $e_0$, $e_8$, $e_5$, $e_2$, $e_{10}$, $e_7$, $e_4$, $e_1$, $e_9$, $e_6$ and $e_3$, and the fold points on $e_0$ have labels $9\p$, $6\p$, $3\p$ and ${11}\p$.
The edge $e_0$ maps $g_{4/11}(e_0)=e_0\bar{e}_8e_8\bar{e}_5e_5\bar{e}_2e_2\bar{e}_{10}e_{10}$.
The labelled tree $G_{4/11}$ is shown in Figure~\ref{fig:star4_11}.
\end{example}

The tree map $g_{1/2}$ is the tree representative for the Smale horseshoe map.
Orbits of $g_{1/2}$ can be coded by taking points in $[v_0,v_2)$ to have code $0$, and points in $(v_2,v_1]$ to have code $1$.
The point $v_2$ has code $\zo$, as it can represent a point in either interval.
We can use the map $h$ in constructing the forcing to give a coding for the map $g_{2/5}$.

In \cite{deCarvalhoHallPPTOPOL}, Section~3, it is shown that the decoration conjecture holds for periodic orbits with star-shaped topological train track.
One of the important stages in the proof is to show that the star map $g_{m/n}$ can be folded and tightened to the star map $g_{p/q}$ if $m/n>p/q$.
However, the methods used are \emph{not} enough to show that the homoclinic braid type $H_{m/n}$ forces $H_{p/q}$.
Using the methods described here, we obtain this result as an immediate consequence of Theorem~\ref{thm:gluetighten}.
The only difficulty is showing that the homoclinic orbit $H_{m/n}$ has topological tree representative $g_{m/n}$.

\begin{theorem}[Tree representative of star homoclinic braids]
Suppose that $0<m/n\leq1/2$, where $m/n$ is a rational expressed in lowest terms.
Then the star homoclinic braid type $H_{m/n}$ with core $c_{m/n}$ has topological tree representative $g_{m/n}$.
\end{theorem}
\begin{proof}
It suffices to embed the Smale horseshoe map, taken as the thickening of the tree map $g_{1/2}$,
 in the thickening of the tree map $g_{m/n}$ and read off the itinerary of the forcing orbit.
An embedding $h_{m/n}$ of $G_{1/2}$ into $g_{m/n}$ is given by taking $G_{1/2}$ to the edge-path 
\[ \alpha_{m/n} e_0\bar{e}_{r}e_{r}\bar{e}_{2r}\cdots \bar{e}_2e_2\cdots e_{1-2r}\bar{e}_{1-r}e_{1-r}\,\bar{e}_{-r}e_{-r}\bar{e}_{-2r}\cdots e_{1+m} \bar{e}_1 , \]
taking indices modulo $n$.
As shown in Section~3.2 of~\cite{deCarvalhoHallPPTOPOL} (with a different labelling of edges),
 the image $g_{m/n}$ of $\alpha_{m/n}$ can be obtained from $\alpha_{m/n}\bar{\alpha}_{m/n}$ by pulling tight.
Hence the thick tree map $\widehat{g}_{m/n}$ of $g_{m/n}$ can be obtained from the Smale horseshoe map $F$ by a pruning isotopy.

Since the only inner fold of $g_{m/n}$ is at the vertex $v_{n\p}$, the map $g_{m/n}$ is the topological tree representative of a trellis
 with a forcing orbit which has the same itinerary as an orbit through $g_{m/n}$.
Hence the core of the forcing orbit $H_{m/n}$ lies in edges $e_1e_0e_{n-1}e_{n-2}\cdots e_3e_2e_1$, so each successive edge is an $m/n$ rotation from the previous edge.
Since the consecutive edges from $e_0$ to $e_2$ give rise to symbol $0$ are those which are $m$ or more edges from $1$, we obtain the word $c_{m/n}$ as required.
Note that in this construction, since the only point of the forcing orbit in $e_2$ is the valence-one vertex $v_2$, we can choose symbol $0$ or $1$; for other orbits, the choice depends on the preimages.
\end{proof}

\begin{theorem}[Forcing for star homoclinic braids]
Suppose that $0<p/q<m/n\leq1/2$, where $p/q$ and $m/n$ are rationals expressed in lowest terms.
Then $H_{m/n}$ forces the star homoclinic braid type $H_{p/q}$.
\end{theorem}

\begin{proof}
The methods of Section~3 of \cite{deCarvalhoHallPPTOPOL} show that $g_{p/q}$ can be obtained from $g_{m/n}$ by glueing and pulling tight.
\end{proof}

Note that the procedure in \cite{deCarvalhoHallPPTOPOL} involves the construction of tree maps $h_{m/n}$ for $0<m/n<1/2$.
While these tree maps are vital to the method of proof, they are \emph{not} the tree represesntatives of any trellis mapping class,
 since there is a fold at the centre vertex.
It should be possible to prove the forcing relation at the trellis level by constructing a sequence of prunings, but a direct trellis construction is not immediately obvious.
Nevertheless, this result shows the power of the trellis theory as an intermediate between (homoclinic) braid types,
 which are the objects of interest, and tree representatives, which are more easily manipulated.

%*****************************************************************************************************************************************************************************

\newcommand{\bibdir}{../bibtex}
\bibliographystyle{alpha}
\bibliography{\bibdir/journalabbreviations,\bibdir/bibliography}

%*****************************************************************************************************************************************************************************

\end{document}

%*****************************************************************************************************************************************************************************